\documentclass[final,3p,times]{elsarticle}

\usepackage{lineno}

\journal{JMAA}

\usepackage{graphicx,amssymb,amsmath,epsf,color,hyperref}
\usepackage{dcolumn}
\usepackage{color}
\usepackage{subfig}
\usepackage{tabularray}
\usepackage{floatrow}
\usepackage{hyperref}

\begin{document}
\title{Efficient and stable derivative-free Steffensen algorithm for root finding}
\date{March 2025}
\author[1]{Alexandre Wagemakers\corref{c1}}
\ead{alexandre.wagemakers@urjc.es}
\affiliation[1]{Nonlinear Dynamics, Chaos and Complex Systems Group, Departamento de Biolog\'ia y Geolog\'ia, F\'isica aplicada y Qu\'imica inorg\'anica, Universidad Rey Juan Carlos, Tulip\'an, M\'ostoles 28933, Madrid, Spain}
\cortext[c1]{Corresponding author}

\author[2]{Vipul Periwal}
\ead{vipulp@mail.nih.gov}
\affiliation[2]{Laboratory of Biological Modeling, National Institute of Diabetes and Digestive and Kidney Diseases, National Institutes of Health, Bethesda, Maryland 20892, USA}

\begin{abstract}
We explore a family of numerical methods, based on the Steffensen divided difference iterative algorithm,  that do not evaluate the derivative of the objective functions. The family of methods achieves second-order convergence with two function evaluations per iteration with marginal additional computational cost. An important side benefit of the method is the improvement in stability for different initial conditions compared to the vanilla Steffensen method. We present numerical results for scalar functions, fields, and scalar fields. This family of methods outperforms the Steffensen method with respect to standard quantitative metrics in most cases.
\end{abstract}

\maketitle

\section{Introduction}

Finding roots of nonlinear functions is a common task in engineering and science. Iterative algorithms are the preferred method to solve such problems, and there is an extensive literature on adapting and implementing methods for specific cases, e.g.~\cite{KhaksarHaghani2015, Piscoran2019, Amat2016, Deng2024}. However, there is still room for improvement in convergence and efficiency for root finding methods. Another challenge is the difficulty of computing or evaluating the derivative of the objective function when it is even possible. In these situations, we must aim for a method avoiding the derivative while keeping the function evaluations as few as possible.

Among the simplest iterative algorithms to find roots without differentiating a function, a popular choice is the Steffensen method. This is a simple modification of the Newton iterative method, where the derivative of the function $f$ is replaced by the first-order divided difference requiring only two function evaluations. The Kung-Traub bound, also known as the Kung-Traub conjecture, states that an optimal iterative method without memory, using $n$ function evaluations per iteration, can achieve a maximum convergence order of $2^{n-1}$~\cite{traub1982iterative}. The Steffensen method is optimal in this sense, achieving an order of 2. This has motivated adaptations of the Steffensen method to multistep schemes to improve the local order of convergence~\cite{cordero2011class, soleymani2011optimal}. The use of memory to store previous values of the evaluations can also improve the convergence of these multipoint methods~\cite{cordero2015efficient, Narang2019,traub1982iterative}.

While the Steffensen method is guaranteed to converge quadratically locally, the estimation of the derivative is far from optimal away from the root. This results in poor stability properties with large regions of initial conditions diverging or converging very slowly. Solutions to this issue have been proposed in~\cite{potra1998secant, amat2006steffensen, hernandez2021algorithm} with modifications to the first-order difference estimation. In particular, these authors noticed that controlling the values of the first-order divided difference improves stability.

Here, we propose a different modification preserving the local convergence behavior while providing a satisfactory estimates of the derivative in regions far from the root. A welcome side-effect of this stability is a reduction in the number of iterations on average. In the present work, we test the method on a collection of nonlinear functions and fields and briefly discuss applications to scalar fields and the implications for optimization. One objective of this article is to show that small modifications to known methods yield significant performance improvements.

\section{A family of Steffensen functions for root finding}
\label{sec:method}

The celebrated Newton-Raphson iteration method is:
\begin{equation}
x_{n+1} = x_n - \frac{f(x_n)}{f'(x_{n})}\ .
\end{equation}
The idea behind the Steffensen-Traub iterative process is to replace the derivative $f'$ with the first-order approximation of the derivative:
\begin{equation}
h(x_n)= \frac{f(x_n + g(f(x_n))) - f(x_n)}{g(f(x_n))}.
\end{equation}
The updated scheme is written as:
\begin{equation}
x_{n+1} = x_n - \frac{f(x_n)}{h(x_{n})}
\end{equation}
When $g(x) = x$ we obtain the Steffensen iterative algorithm requiring two evaluations of the function $f$ at each step. While the method guarantees second-order convergence in the vicinity of a root, the estimation of the derivative is rather poor for initial conditions far from the solution. These inaccurate estimates lead to stability problems and divergent trajectories. 

In this article, The proposed functions $g$ solve the stability problem and improve the speed of convergence. A direct consequence of our modification is a reduction in the number of iteration and an increased range of valid initial conditions. The local quadratic convergence for smooth functions still holds for a large class of functions $g$. We will show empirically that the choice of nonlinear functions performs better in numerical simulations. In the articles~\cite{amat2006steffensen, amat2015improving}, a similar argument has been proposed using a basic control algorithm to keep the value of $g(f(x_n))$ within bounds. We generalize their approach baking the ad hoc argument into the function $g$.  

An interesting modification of the Steffensen algorithm has been introduced in several works with a constant $\gamma$ in the estimation of the derivative: 
\begin{equation}
   \Gamma_n = \frac{f(x_n +\gamma f(x_n)) - f(x_n)}{\gamma f(x_n)}.
\end{equation}
The properties of the iteration scheme are well studied, see for example \cite{traub1982iterative, cordero2015efficient}. If we treat the constant $\gamma$ as a dynamical variable depending on the previous iteration we can obtain an accelerated Steffensen method with convergence order at least $1 + \sqrt{2}$. This has been described in~\cite{traub1982iterative} section 8.6.  

The accelerated algorithm is summarized as:
\begin{align}
    &x_{n+1} = x_n - \frac{f(x_n)}{\Gamma_n},\\
    &\Gamma_n = \frac{f(x_n + g(\gamma_n f(x_n))) - f(x_n)}{g(\gamma_n f(x_n))}, \\
    &\gamma_n = -\frac{1}{\Gamma_{n-1}}.
\end{align}
The difference with the normal method is an additional cost in memory storage to keep the previous value of $\Gamma_k$ at each step and two multiplications.These computation costs actually increase quadratically with the dimension of the domain of the function. Notice that the accelerated method can be integrated in a generalized function $\varphi(u,v) = g(u\cdot v)$. This trivial definition is useful for the generalization of the computer code abstraction. 

Considering these two numerical schemes, we will compare several functions $g$ for the approximation of the derivative, Jacobians, and gradients to obtain better convergence and stability results.

\section{\label{sec:order}Order of convergence}

The convergence of the scheme is of second order for smooth functions. We can prove this convergence by approximating the error $e_{k} = x_k - x^*$ up to the second order.

We define:
\begin{equation}
    v(x) = x - \frac{f(x)}{h(x)},
\end{equation}
with $h(x)$ our derivative estimate. We expand the function $v$ using a Taylor expansion near a root $x^*$:
\begin{equation}
    v(x_k) = v(e_k + x^*) = v(x^*) + v'(x^*) e_k + \frac{1}{2} v''(x^*) e_k^2 + O(e_k^3),
\end{equation}
with $v(x^*) = x^*$. To prove the convergence, it is sufficient to prove that $v'(x^*) = 0$ and $v''(x^*) \neq 0$. For the first part, we differentiate:
\begin{equation}
    v'(x) = 1 - \frac{f'(x) h(x) - h'(x)f(x)}{h^2(x)}.
\end{equation}
The evaluation near the root gives the expression:
\begin{equation}
    v'(x^*) =  1 - \frac{g(f(x^*))f'(x^*)}{f(x^* + g(f(x^*))}.
\end{equation}
For a smooth function $g$ such that $g(0)= 0$ and $g'(0) \neq 0$, we can expand $f$ again around $g(f(x^*))$ at first order as $f(x^* + g(f(x^*))) \simeq f(x^*) + g(f(x^*)) f'(x^*)$, obtaining:
\begin{equation}
    v'(x^*) =   1 - \frac{g(f(x^*)) f'(x^*)}{g(f(x^*)))f'(x^*)} = 0.
\end{equation}
If the estimated derivative $h$ converges to the true derivative $f'$, the previous expression is meaningful. The last step consists in writing $v(x_n) - x^* = e_{k+1}$ and we obtain the recurrence formula for the error:
\begin{equation}
    e_{k+1} = \frac{1}{2}v''(x^*) e_k^2 + O(e_k^3)
\end{equation}
The evaluation of $v''(x^*)= f''(x^*)/f'(x^*) \not= 0$ follows from the simple root assumption and 
\begin{equation}
    v''(x) = - \frac{f''(x) h(x) - h''(x)f(x)}{h^2(x)} + 2 \frac{h(x)}{h'(x)} \left(1-v'(x)\right)
\end{equation}
evaluated at $x = x^*.$
This proves the convergence for smooth functions. A more general result proven in~\cite{candela2019class} shows that any iteration scheme in the form $x_{n+1} = x_n - f(x_n)/\hat f'(x_n)$ where $\hat f'$ is an approximation to $f'$ has convergence order $\min((x - x^*)(f'(x^*) - \hat f'(x^*)), 2)$. The theorem holds under reasonable assumptions on the regularity of $f$. Another argument for the convergence of semi-smooth functions has been developed in~\cite{potra1998secant} and~\cite{amat2006steffensen} and can be adapted to our scheme.

The order of convergence of the accelerated method has been demonstrated in \cite{traub1982iterative,Zhao2024} and will not be demonstrated here. We will rely on the fact that the selected functions $g$ are well approximated by the identity function near the origin. Our purpose is to accelerate the convergence using a better approximation of the derivative for the first steps of the algorithm.

To confirm numerically the quadratic convergence near the origin, we estimate the order near the root with the ACOC convergence estimator defined as \cite{grau2010}:
\begin{equation}
q_n = \frac{\log(|x_{n+1}-x_n|/|x_n - x_{n-1}|)}{\log(|x_n-x_{n-1}|/|x_{n-1} - x_{n-2}|)}.
\end{equation}
The estimation stops until the criterion $|x_{n+1} - x_n|<\varepsilon$ for $\varepsilon = 10^{-25}$ has been reached and we pick the last value of the series $q_n$. We present the results in the numerical section.

\section{\label{sec:method}Solving roots for functions in $\mathbb{R}^k \to \mathbb{R}^k$}

Newton's classical method in higher dimensions~\cite{ortega2000iterative} is expressed as:
\begin{equation}
     \bar x_{n+1} = \bar  x_n - J_f(\bar x_n)^{-1} f(\bar x_n)
\end{equation}
where $J_f$ is the Jacobian matrix of $f$. We define the operator:
\begin{equation}
    T_{i,j}(\bar x) = \frac{f_i(x_1, \dots , x_j + g(f_i(\bar x)), \dots, x_k) - f_i(\bar x)}{g(f_i(\bar x))}
\end{equation}

We can extend our method with an estimate of the Jacobian matrix such that the partial derivatives are estimated with the operator $T$:
\begin{equation}
    \left. \frac{\partial f_i}{\partial x_j} \right |_{\bar x}  \simeq T_{i,j}(\bar x)
\end{equation}
where $\bar x$ is the coordinate vector in $\mathbb{R}^k$. The Jacobian matrix can be evaluated with $k^2$ extra evaluations of the components $f_i$:
\begin{equation}
    J_f(\bar x)  \equiv
\begin{bmatrix}
    T_{1,1}(\bar x) &     T_{1,2}(\bar x) &  \dots  & T_{1,k}(\bar x)   \\

    T_{2,1}(\bar x) &     T_{2,2}(\bar x) &  \dots  & T_{2,k}(\bar x)   \\

    \vdots & \vdots & \ddots & \vdots \\
    T_{k,1}(\bar x) &     T_{k,2}(\bar x) &  \dots  & T_{k,k}(\bar x)
\end{bmatrix}
\end{equation}

At each step, the algorithm requires $k^2 + k$ function evaluations. 

This reasoning can be extended to the estimation of gradients of scalar functions from $\mathbb{R}^k \to \mathbb{R}$. Using the same notation, we have the gradient:
\begin{equation}\label{eq:grad_est}
    \nabla f(\bar x) \simeq (T_{1,1}(\bar x), T_{1,2}(\bar x), \cdots ,  T_{1,k}(\bar x))
\end{equation}
This approximation is useful in the search for roots of these functions with an adaptation of the Steffensen algorithm.

\section{\label{sec:num}Numerical results}

To compare methods, we use a collection of scalar functions and fields listed in Tab.~\ref{tab1}, mostly adapted from~\cite{cordero2007variants}. 

\begin{table}[]
    \centering
\begin{tabular}{|l|}
\hline
$f_1(x) = x^3 - 9x^2 + 28x - 30$ \\
$f_2(x) = \sin(x) + x \cos(x)$ \\
$f_3(x) = e^{x^2} - e^{\sqrt{2} x}$ \\
$f_4(x) = (\sin(x) - x/2)^2$ \\
$f_5(x) = \tan^{-1}(x)$ \\
$f_6(x) = (x - 1)^6 - 1$ \\
$f_7(x) = 4 \sin(x) - x + 1$ \\
$f_8(x) = (x^2 - 1)(x^2 + 1)$ \\
$f_9(x) = (x^2 - 4)(x + 1.5)(x - 0.5)$ \\
$f_{10}(x) = (x + 2)(x + 1.5)^2(x - 0.5)(x - 2)$ \\
$f_{11}(x) = (x - 1)^3 + 4(x - 1)^2 - 10$ \\
$f_{12}(x) = \sin\left(x - \frac{14}{10}\right)^2 - \left(x - \frac{14}{10}\right)^2 + 1$ \\
$f_{13}(x) = x^2 - e^x - 3x + 2$ \\
$f_{14}(x) = \left(x + \frac{5}{4}\right) e^{\left(\left(x + \frac{5}{4}\right)^2\right)} - \sin\left(x + \frac{5}{4}\right)^2 + 3 \cos\left(x + \frac{5}{4}\right) + 5$ \\
$f_{15}(x, y) = \left( x + e^{y} - \cos(y), 3x - y - \sin(y) \right)$ \\
$f_{16}(x, y) = \left( e^{x^2} + 8x \sin(y), x + y - 1 \right)$ \\
$f_{17}(x, y) = \left( \sin(x) + y \cos(x), x - y \right)$ \\
$f_{18}(x, y) = \left( x^2 - 2x - y + 0.5, x^2 + 4y^2 - 4.0 \right)$ \\
$f_{19}(x, y) = \left( e^{x^2} - e^{\sqrt{2} x}, x - y \right)$ \\
$f_{20}(x, y, z, w) = \left( y z + w (y + z), x z + w (x + z),\right.$ \\
~~~~~$\left. x y + w (x + y), x y + x z + y z - 1 \right)$ \\
$f_{21}(x, y) = \left( \Re\left\{(x + i y)^3 - 1\right\}, \Im\left\{(x + i y)^3 - 1\right\} \right)$\\
\hline
\end{tabular}
\caption{List of functions tested}
\label{tab1}
\end{table}

First, we evaluate  the methods using a single initial condition using the multi-precision floating-point arithmetic library (GNU MPFR). This library is
integrated in the Julia programming language base package. All the numerical computations have been performed in Julia using standard packages and the scripts available at: \url{https://github.com/awage/SteffensenRootFinding}.

The stopping criterion is $|f(x_n)| < \varepsilon$ for $\varepsilon = 1\cdot 10^{-25}$ and we use floating-point numbers up to 100 decimal digits. Table~\ref{tab2} shows the number of iterations and the final root for a single initial condition for the following functions $g_i$:
\begin{itemize}
    \item $g_1(x) = \tanh(x)$
    \item $g_2(x) = \max(|x|, 1) \times \text{sign}(x)$
\end{itemize}
The function $g_1$ and $g_2$ are similar in shape but $g_2$ is significantly faster to compute. Nevertheless, results are often better with the continuous and differentiable function $g_1$.

\begin{table*}[]
    \centering
    \begin{tabular}{p{0.5cm}p{1cm} lll | lll | lll |}
 & & \multicolumn{3}{c}{It.} & \multicolumn{3}{c}{Root}  & \multicolumn{3}{c}{q}  \\
 & & $g_1$ & $g_2$ & SM &
    $g_1$ & $g_2$ & SM &
    $g_1$ & $g_2$ & SM\\
 \hline
{\footnotesize $f_{1}$}
& {\footnotesize (norm.)}
    & {\bf 11} & 12 & 22 
 & 3.0  & 3.0  & 3.0  
 & 3.2 & 3.0 & 3.0 \\
& {\footnotesize (accel.)}
    & {\bf 7} & {\bf 7} & 8 
 & 3.0  & 3.0  & 3.0  
 & 3.7 & 3.7 & 3.7 \\
{\footnotesize $f_{2}$}
& {\footnotesize (norm.)}
    & {\bf 6} & {\bf 6} & 8 
 & 2.03  & 2.03  & 2.03  
 & 2.0 & 2.0 & 2.0 \\
& {\footnotesize (accel.)}
    & 7 & {\bf 6} & {\bf 6} 
 & 2.03  & 2.03  & 4.91  
 & 3.6 & 1.3 & 3.1 \\
{\footnotesize $f_{3}$}
& {\footnotesize (norm.)}
    & {\bf 17} & 23 & 25 
 & 1.41  & 1.41  & 1.41  
 & 2.0 & 2.0 & 2.0 \\
& {\footnotesize (accel.)}
    & {\bf 8} & 16 & 9 
 & 1.41  & 1.41  & 1.41  
 & 3.0 & 2.4 & 4.6 \\
{\footnotesize $f_{4}$}
& {\footnotesize (norm.)}
 & 43 & 43 & 43 
 & -1.9  & -1.9  & -1.9  
 & 1.0 & 1.0 & 1.0 \\
& {\footnotesize (accel.)}
 & 34 & 34 & 34 
 & 0.0  & 0.0  & 0.0  
 & 1.0 & 1.0 & 1.0 \\
{\footnotesize $f_{5}$}
& {\footnotesize (norm.)}
 & 6 & 6 & 6 
 & 0.0  & 0.0  & 0.0  
 & 2.9 & 2.9 & 2.9 \\
& {\footnotesize (accel.)}
 & 4 & 4 & 4 
 & 0.0  & 0.0  & 0.0  
 & - & - & - \\
{\footnotesize $f_{6}$}
& {\footnotesize (norm.)}
    & {\bf 48} & 52 & nc 
 & 2.0  & 2.0  & nc  
 & 2.0 & 2.0 & 0.0 \\
& {\footnotesize (accel.)}
    & {\bf 8} & 12 & 26 
 & 2.0  & 2.0  & 2.0  
 & 1.4 & 1.5 & 2.4 \\
{\footnotesize $f_{7}$}
& {\footnotesize (norm.)}
 & 6 & 6 & 6 
 & 2.7  & 2.7  & 2.7  
 & 2.0 & 2.0 & 2.0 \\
& {\footnotesize (accel.)}
 & 6 & 6 & 6 
 & 2.7  & 2.7  & 2.7  
 & 2.8 & 3.1 & 3.1 \\
{\footnotesize $f_{8}$}
& {\footnotesize (norm.)}
 & {\bf 16} & 17 & 152 
 & 1.0  & 1.0  & 1.0  
 & 2.0 & 2.0 & 2.0 \\
& {\footnotesize (accel.)}
 & {\bf 6} & {\bf 6} & 9 
 & 1.0  & 1.0  & 1.0  
 & 2.6 & 2.6 & 2.4 \\
{\footnotesize $f_{9}$}
& {\footnotesize (norm.)}
 & {\bf 7} & {\bf 7} & nc 
 & 0.5  & 0.5  & nc
 & 2.0 & 2.0 & - \\
& {\footnotesize (accel.)}
 & 7 & 7 & 7 
 & 2.0  & 2.0  & 0.5  
 & 2.0 & 2.1 & 2.6 \\
{\footnotesize $f_{10}$}
& {\footnotesize (norm.)}
 & {\bf 31} & 32 & nc 
 & 2.0  & 2.0  & nc
 & 2.0 & 2.0 & 0.0 \\
& {\footnotesize (accel.)}
 & {\bf 12} & {\bf 12} & 13 
 & 2.0  & 2.0  & 2.0  
 & 2.4 & 2.4 & 2.4 \\
{\footnotesize $f_{11}$}
& {\footnotesize (norm.)}
 & {\bf 12} & 13 & 242 
 & 2.37  & 2.37  & 2.37  
 & 2.0 & 2.0 & 2.0 \\
& {\footnotesize (accel.)}
 & {\bf 7} & {\bf 7} & 9 
 & 2.37  & 2.37  & 2.37  
 & 2.4 & 2.4 & 2.4 \\
{\footnotesize $f_{12}$}
& {\footnotesize (norm.)}
 & {\bf 8} & 9 & nc 
 & 2.8  & 2.8  & nc  
 & 2.0 & 2.0 & - \\
& {\footnotesize (accel.)}
 & {\bf 7} & {\bf 7} & 8 
 & 2.8  & 2.8  & 2.8  
 & 2.3 & 2.3 & 2.4 \\
{\footnotesize $f_{13}$}
& {\footnotesize (norm.)}
 & {\bf 7} & {\bf 7} & 8 
 & 0.26  & 0.26  & 0.26  
 & 2.0 & 2.0 & 2.0 \\
& {\footnotesize (accel.)}
 & {\bf 6} & {\bf 6} & 7 
 & 0.26  & 0.26  & 0.26  
 & 2.8 & 2.8 & 2.3 \\
{\footnotesize $f_{14}$}
& {\footnotesize (norm.)}
 & nc & nc & nc 
 & nc  & nc  & nc
 & - & - & - \\
& {\footnotesize (accel.)}
 & 10 & {\bf 9} & nc 
 & -2.46  & -2.46  & 0.3  
 & 2.4 & 2.4 & - \\
{\footnotesize $f_{15}$}
& {\footnotesize (norm.)}
 & {\bf 10} & {\bf 10} & 19 
 & (0.0, 0.0) & (0.0, 0.0) & (0.0, 0.0) 
 & 2.0 & 2.0 & 2.0 \\
& {\footnotesize (accel.)}
 & {\bf 7} & {\bf 7} & nc 
 & (0.0, 0.0, ) & (0.0, 0.0, ) & nc
 & 1.9 & 1.9 & - \\
{\footnotesize $f_{16}$}
& {\footnotesize (norm.)}
 & {\bf 7} & {\bf 7} & nc 
 & (-0.14, 1.14) & (-0.14, 1.14) & nc
 & 2.0 & 2.0 & - \\
& {\footnotesize (accel.)}
 & {\bf 8} & {\bf 8} & nc 
 & (-0.14, 1.14) & (-0.14, 1.14) & nc
 & 2.0 & 2.0 & - \\
{\footnotesize $f_{17}$}
& {\footnotesize (norm.)}
 & {\bf 5} & 6 & nc 
 & (0.0, 0.0) & (0.0, 0.0) & nc
 & 2.0 & 3.0 & - \\
& {\footnotesize (accel.)}
 & {\bf 4} & 5 & nc 
 & (0.0, 0.0) & (0.0, 0.0) & nc
 & - & 3.0 & - \\
{\footnotesize $f_{18}$}
& {\footnotesize (norm.)}
 & {\bf 17} & nc & nc 
 & (-0.22, 0.99) & nc & nc
 & 2.0 & - & - \\
& {\footnotesize (accel.)}
 & 10 & {\bf 8} & {\bf 8} 
 & (-0.22, 0.99) & (-0.22, 0.99) & (-0.22, 0.99) 
 & 1.7 & 1.9 & 1.9 \\
{\footnotesize $f_{19}$}
& {\footnotesize norm.}
 & nc & nc & nc 
 & (7.0, 7.0) & (7.0, 7.0) & (7.0, 7.0) 
 & 0.0 & 0.0 & 0.0 \\
& {\footnotesize accel.}
 & {\bf 37} & {\bf 37} & nc 
 & (1.41, 1.41) & (1.41, 1.41) & (7.0, 7.0) 
 & 2.4 & 2.4 & 0.0 \\
{\footnotesize $f_{20}$}
& {\footnotesize (norm.)}
 & 7 & 7 & 7 
 & $\alpha$ & $\alpha$ & $\alpha$ 
 & 2.2 & 2.2 & 2.2 \\
& {\footnotesize (accel.)}
 & 7 & 7 & 7 
 & $\alpha$ & $\alpha$ & $\alpha$ 
 & 2.2 & 2.2 & 2.2 \\
{\footnotesize $f_{21}$}
& {\footnotesize (norm.)}
 & 14 & {\bf 13} & nc
 & (-0.5, 0.87) & (-0.5, 0.87) & nc
 & 2.0 & 2.0 & - \\
& {\footnotesize (accel.)}
 & {\bf 10} & {\bf 10} & 19 
 & (1.0, -0.0) & (1.0, -0.0) & (1.0, 0.0) 
 & 1.8 & 1.8 & 1.7 \\
\end{tabular}

    \caption{\label{tab2}Number of iterations and final root for a single initial condition. The Steffensen Method (SM) is tested against the modified version with the functions $g_i$. Non-convergent initial conditions are noted as nc and the evaluation of the convergence order requires at least 6 iterations. The root $\alpha$ in the table is $\alpha = (0.58,0.58, 0.58, -0.29)$. Values are represented up to the second decimal.}
\end{table*}

We introduce a small modification in $g_1$ and $g_2$ to enforce numerical stability when the function $f$ is very close to a root. The estimation of the derivative can fail due to the limit of available digits. A simple modification consists of setting a lower threshold such that if $|x| < \delta$, $g$ return $\pm \delta$. For the simulations, we have set $\delta = \varepsilon/2$. The algorithm is iterated for a maximum of 1000 iterations, and beyond this value, we mark the initial condition as non-convergent. Finally, the table includes the estimation of the convergence order $q$ computed with the ACOC method detailed in Sec.~\ref{sec:order}. For each function, we compute the normal and accelerated Steffensen method, marked as SM in the table, with and without the functions $g_i$.

From the values in Tab.~\ref{tab2}, a clear result appears at first sight. The modified methods are better in terms of the number of iterations and stability. The estimated order of convergence is also consistent with the predicted order of convergence, although there are exceptions. In particular, function $f_4$ shows linear convergence due to the presence of multiple roots. The results for the two functions are similar, although there is a slight edge for the function $g_1$.

\begin{table}[]
    \centering
    \begin{tabular}{|p{0.4cm}p{1cm} lll | lll | lll |}
    \hline
 & & \multicolumn{3}{c}{Non conv.} & \multicolumn{3}{|c}{Mean it.}  & \multicolumn{3}{|c|}{T/it.}  \\
 & & $g_1$ & $g_2$ & SM &
    $g_1$ & $g_2$ & SM &
    $g_1$ & $g_2$ & SM\\
 \hline
{\footnotesize $f_{1}$}
& {\footnotesize (norm.)}
 & \textbf{0.0} & \textbf{0.0} & 51.1 
 & 12.3 & 12.5 & 52.2 
 & 1.61 & \textbf{0.25} & 1.0 \\
& {\footnotesize (accel.)}
 & \textbf{0.0} & \textbf{0.0} & \textbf{0.0} 
    & \textbf{7.4} & \textbf{7.4} & 8.4 
    & 0.43 & \textbf{0.41} & 1.0 \\
    \hline
{\footnotesize $f_{2}$}
& {\footnotesize (norm.)}
 & \textbf{0.0} & \textbf{0.0} & \textbf{0.0} 
    & \textbf{5.4} & 5.5 & 6.2 
    & \textbf{0.99} & 3.58 & 1.0 \\
& {\footnotesize (accel.)}
 & \textbf{0.0} & \textbf{0.0} & \textbf{0.0} 
    & \textbf{4.6} & \textbf{4.6} & 5.1 
    & 0.53 & \textbf{0.47} & 1.0 \\
    \hline
{\footnotesize $f_{3}$}
& {\footnotesize (norm.)}
 & 39.0 & 39.2 & \textbf{12.3} 
 & 14.2 & 14.2 & \textbf{2.4} 
 & \textbf{0.37} & \textbf{0.37} & 1.0 \\
& {\footnotesize (accel.)}
 & 0.3 & \textbf{0.2} & 3.9 
    & 30.6 & 30.3 & \textbf{3.3} 
    & 3.79 & 3.4 & \textbf{1.0} \\
    \hline
{\footnotesize $f_{4}$}
& {\footnotesize (norm.)}
    & 8.9 & \textbf{8.8} & 9.1 
    & 26.4 & 26.3 & \textbf{17.0} 
    & 1.32 & 1.29 & \textbf{1.0} \\
& {\footnotesize (accel.)}
 & 5.5 & 6.1 & \textbf{0.0} 
    & 20.3 & 21.6 & \textbf{15.5} 
    & 0.58 & \textbf{0.55} & 1.0 \\
    \hline
{\footnotesize $f_{5}$}
& {\footnotesize (norm.)}
    & \textbf{91.0} & 91.6 & 91.6 
    & 5.0 & \textbf{4.9} & \textbf{4.9} 
    & 1.26 & 1.09 & \textbf{1.0} \\
& {\footnotesize (accel.)}
 & 78.9 & 77.7 & \textbf{0.0} 
    & 4.2 & \textbf{4.1} & 6.7 
    & \textbf{0.68} & 2.54 & 1.0 \\
    \hline
{\footnotesize $f_{6}$}
& {\footnotesize (norm.)}
 & \textbf{0.0} & \textbf{0.0} & 91.7 
    & 25.9 & 27.5 & \textbf{19.0} 
    & \textbf{0.63} & 0.67 & 1.0 \\
& {\footnotesize (accel.)}
 & \textbf{0.0} & \textbf{0.0} & 59.3 
    & 11.9 & 11.7 & \textbf{9.6} 
    & 0.52 & \textbf{0.44} & 1.0 \\
    \hline
{\footnotesize $f_{7}$}
& {\footnotesize (norm.)}
 & 22.8 & 22.8 & \textbf{0.0} 
 & 9.9 & 10.1 & \textbf{7.2} 
 & 1.28 & 1.34 & \textbf{1.0} \\
& {\footnotesize (accel.)}
 & 28.6 & 30.0 & \textbf{0.0} 
 & 13.7 & 14.6 & \textbf{6.6} 
 & \textbf{0.41} & 0.43 & 1.0 \\
    \hline
{\footnotesize $f_{8}$}
& {\footnotesize (norm.)}
 & \textbf{0.0} & \textbf{0.0} & 86.0 
 & \textbf{13.9} & 14.2 & 26.4 
 & 0.46 & \textbf{0.42} & 1.0 \\
& {\footnotesize (accel.)}
 & \textbf{0.0} & \textbf{0.0} & 0.5 
 & \textbf{8.8} & \textbf{8.8} & 9.8 
 & 0.41 & \textbf{0.4} & 1.0 \\
    \hline
{\footnotesize $f_{9}$}
& {\footnotesize (norm.)}
 & \textbf{0.0} & \textbf{0.0} & 82.1 
 & \textbf{12.5} & 13.1 & 19.6 
 & 0.5 & \textbf{0.48} & 1.0 \\
& {\footnotesize (accel.)}
 & \textbf{0.0} & \textbf{0.0} & 0.3 
 & \textbf{8.0} & \textbf{8.0} & 9.0 
 & 0.43 & \textbf{0.39} & 1.0 \\
    \hline
{\footnotesize $f_{10}$}
& {\footnotesize (norm.)}
 & \textbf{0.0} & \textbf{0.0} & 85.4 
 & \textbf{25.3} & 26.1 & 40.1 
 & 0.23 & \textbf{0.22} & 1.0 \\
& {\footnotesize (accel.)}
 & \textbf{0.0} & \textbf{0.0} & 2.4 
 & 9.5 & \textbf{9.4} & 10.6 
 & \textbf{0.39} & \textbf{0.39} & 1.0 \\
    \hline
{\footnotesize $f_{11}$}
& {\footnotesize (norm.)}
 & 15.4 & \textbf{5.4} & 53.3 
 & \textbf{28.8} & 45.9 & 61.5 
 & 0.68 & \textbf{0.54} & 1.0 \\
& {\footnotesize (accel.)}
 & 0.1 & \textbf{0.1} & 2.5 
 & 23.6 & \textbf{22.2} & 35.5 
 & 0.22 & \textbf{0.18} & 1.0 \\
    \hline
{\footnotesize $f_{12}$}
& {\footnotesize (norm.)}
 & \textbf{0.0} & \textbf{0.0} & 34.3 
 & \textbf{7.8} & 7.9 & 13.1 
 & 0.29 & \textbf{0.27} & 1.0 \\
& {\footnotesize (accel.)}
 & \textbf{0.0} & \textbf{0.0} & \textbf{0.0} 
 & \textbf{6.0} & \textbf{6.0} & 7.2 
 & 0.49 & \textbf{0.42} & 1.0 \\
    \hline
{\footnotesize $f_{13}$}
& {\footnotesize (norm.)}
 & \textbf{0.0} & \textbf{0.0} & 38.1 
 & \textbf{6.2} & 6.4 & 11.5 
 & \textbf{0.67} & 0.68 & 1.0 \\
& {\footnotesize (accel.)}
 & \textbf{0.0} & \textbf{0.0} & 14.1 
 & 6.7 & 6.6 & \textbf{5.5} 
 & \textbf{0.67} & 2.38 & 1.0 \\
    \hline
{\footnotesize $f_{14}$}
& {\footnotesize (norm.)}
 & 90.4 & 90.6 & 99.2 
 & 85.4 & 89.4 & \textbf{26.1} 
 & \textbf{0.78} & 0.81 & 1.0 \\
& {\footnotesize (accel.)}
 & \textbf{4.0} & \textbf{4.0} & 13.0 
 & 35.1 & 34.9 & \textbf{3.1} 
 & 4.4 & 4.2 & \textbf{1.0} \\
    \hline
    \hline
{\footnotesize $f_{15}$}
& {\footnotesize (norm.)}
 & \textbf{11.6} & \textbf{11.6} & 74.3 
 & 33.9 & \textbf{33.7} & 38.9 
 & \textbf{0.89} & 0.93 & 1.0 \\
& {\footnotesize (accel.)}
 & \textbf{8.5} & \textbf{8.5} & 82.1 
 & 29.9 & 29.8 & \textbf{11.3} 
 & 2.81 & 2.52 & \textbf{1.0} \\
    \hline
{\footnotesize $f_{16}$}
& {\footnotesize (norm.)}
 & \textbf{42.4} & \textbf{42.4} & 99.9 
 & 28.5 & 32.0 & \textbf{6.7} 
 & 1.36 & 1.62 & \textbf{1.0} \\
& {\footnotesize (accel.)}
 & \textbf{12.3} & 16.3 & 99.9 
 & 46.3 & 44.6 & \textbf{6.7} 
 & 2.43 & 2.4 & \textbf{1.0} \\
    \hline
{\footnotesize $f_{17}$}
& {\footnotesize (norm.)}
 & \textbf{0.0} & \textbf{0.0} & 99.9 
 & 6.1 & 6.1 & \textbf{3.0} 
 & 1.63 & 1.53 & \textbf{1.0} \\
& {\footnotesize (accel.)}
 & \textbf{0.0} & \textbf{0.0} & 98.0 
 & \textbf{5.7} & 5.8 & 12.0 
 & 0.94 & \textbf{0.44} & 1.0 \\
    \hline
{\footnotesize $f_{18}$}
& {\footnotesize (norm.)}
 & \textbf{0.0} & \textbf{0.0} & 68.3 
 & 15.3 & \textbf{14.5} & 31.5 
 & \textbf{0.32} & 0.34 & 1.0 \\
& {\footnotesize (accel.)}
 & \textbf{0.0} & \textbf{0.0} & \textbf{0.0} 
 & 17.3 & 18.4 & \textbf{17.1} 
 & 1.21 & 1.21 & \textbf{1.0} \\
    \hline
{\footnotesize $f_{19}$}
& {\footnotesize (norm.)}
 & \textbf{38.9} & 39.0 & 95.7 
 & 14.1 & 14.2 & \textbf{3.0} 
 & 1.82 & 1.06 & \textbf{1.0} \\
& {\footnotesize (accel.)}
 & 0.5 & \textbf{0.3} & 97.9 
 & 30.6 & 30.3 & \textbf{3.5} 
 & 2.13 & 1.01 & \textbf{1.0} \\
    \hline
{\footnotesize $f_{20}$}
& {\footnotesize (norm.)}
 & \textbf{54.9} & 55.1 & 55.1 
 & 19.7 & \textbf{19.3} & 19.4 
 & 1.03 & 1.02 & \textbf{1.0} \\
& {\footnotesize (accel.)}
 & 55.2 & 55.0 & \textbf{55.1} 
 & \textbf{19.2} & 19.4 & 19.4 
 & 1.12 & \textbf{1.0} & \textbf{1.0} \\
    \hline
{\footnotesize $f_{21}$}
& {\footnotesize (norm.)}
 & \textbf{0.0} & \textbf{0.0} & 98.2 
 & \textbf{15.3} & 15.6 & 18.6 
 & 0.49 & \textbf{0.48} & 1.0 \\
& {\footnotesize (accel.)}
 & \textbf{0.0} & \textbf{0.0} & \textbf{0.0} 
 & \textbf{11.7} & \textbf{11.7} & 13.9 
 & 0.87 & \textbf{0.79} & 1.0 \\
 \hline
\end{tabular}

    \caption{\label{tab3}Percentage of non-converging initial conditions, average number of iterations, and compute time per iteration for $10^4$ sampled initial conditions in the interval $[-10, 10]^d$ where $d$ is the dimension of the function.}
\end{table}

This local picture of the algorithm reflects only partially the properties of the methods. From the perspective of applications, crucial differences appear while examining the convergence and number of iterations over a larger domain.  

Table~\ref{tab3} describes averages on initial conditions sampled in the region $[-10, 10]^d$ where $d$ is the dimension of the domain of the nonlinear function. Three different figures of merit are computed: the percentage of non-convergent initial conditions, the average number of iterations for convergent initial conditions, and finally the computation time per iteration normalized to the Steffensen algorithm value. For this last metric, we have compared the values only when the three methods converged to a root for a given initial condition. For this simulation, we have used standard double-precision arithmetic with the stopping criterion set to $10^{-8}$ to benchmark the use case of a practical application of the algorithm. The maximum number of iterations is set to 200. 

The results show a notable increase in stability of the algorithm in terms of convergence for the modified method. The number of iterations needed is in most cases lower. The computation time per iteration depends on several factors: the number of iterations and the computation of the function $g$. The evaluation of $g_1$ is expensive, in particular, the hyperbolic tangent depends on many arithmetic operations. Even with this limitation, the numbers are in favor of the function $g_1$ in many cases. The most favorable results have been highlighted in bold face in the table.

\section{Dynamics on the plane}

The iterated algorithm can also be studied from the dynamical system point of view with tools from nonlinear dynamics such as the basins of attraction~\cite{jo2024annealing}. For applications from the plane to the plane, a common representation is to associate an initial condition to a root forming a basin of attraction. Structures in the phase plane such as the boundaries appear separating basins.

Basins for maps from $\mathbb{R}^2$ to $\mathbb{R}^2$ can be directly represented graphically. We can also adapt a function $f$ to the complex plane defining the two-dimensional function $F(x,y) = (\Re\{f(x+iy)\} , \Im\{f(x+iy)\})$. This transformation allows the estimation of the Jacobian matrix with the presented method in Sec.\ref{sec:num} for the new function $F(x) = (h(x,y), j(x,y))$, where $h$ represents the real part of $f(z)$ and $j$ its imaginary part.

\begin{figure*}
\begin{tblr}{
  colspec = {X[c]X[c,h]X[c,h]X[c,h]},
  stretch = 0,
  rowsep = 2pt,
}
    & $g_1$ & $g_2$ & SM\\
$f_9$ normal  &\includegraphics[width = 0.2\textwidth]{"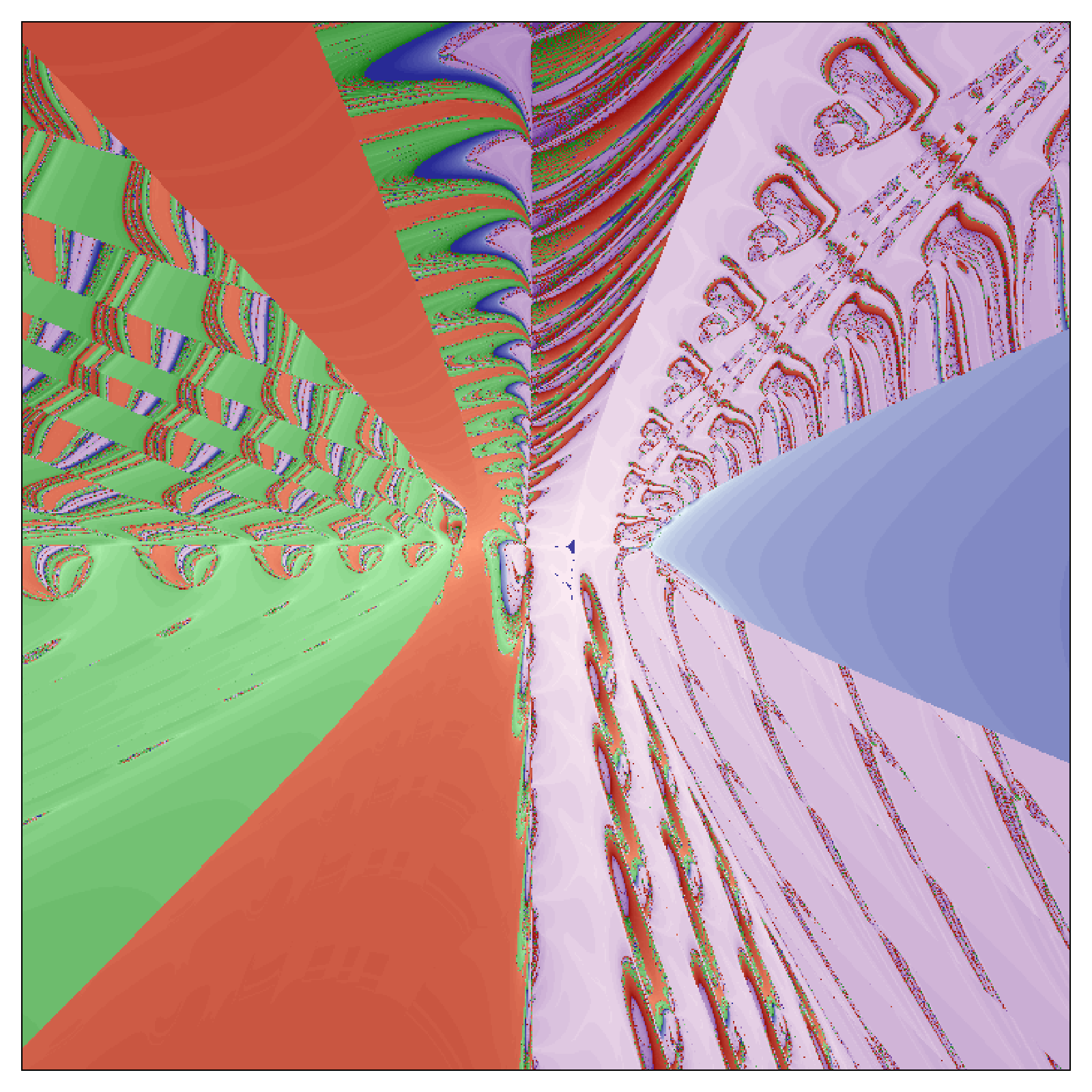"} & 
    \includegraphics[width = 0.2\textwidth]{"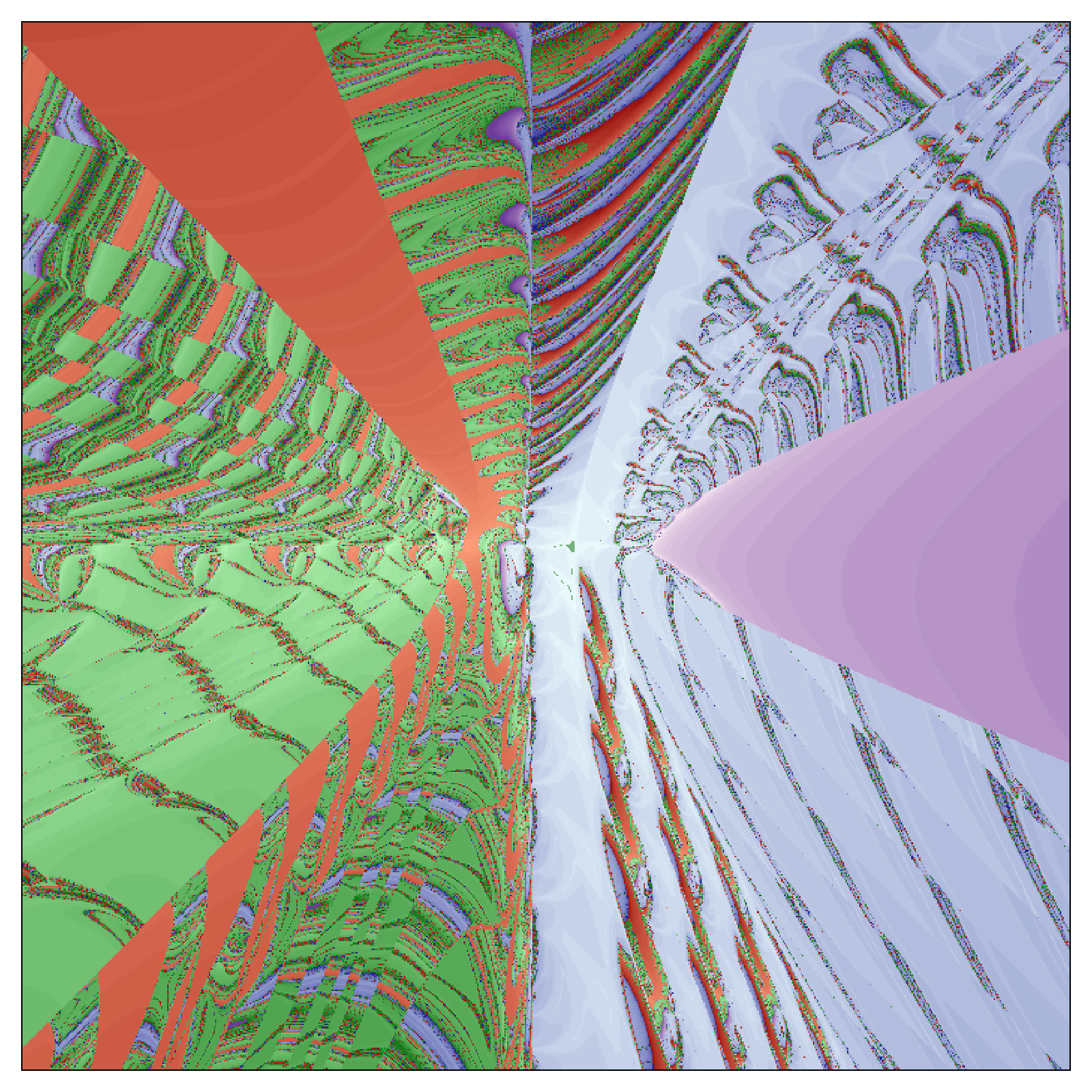"} &
    \includegraphics[width = 0.2\textwidth]{"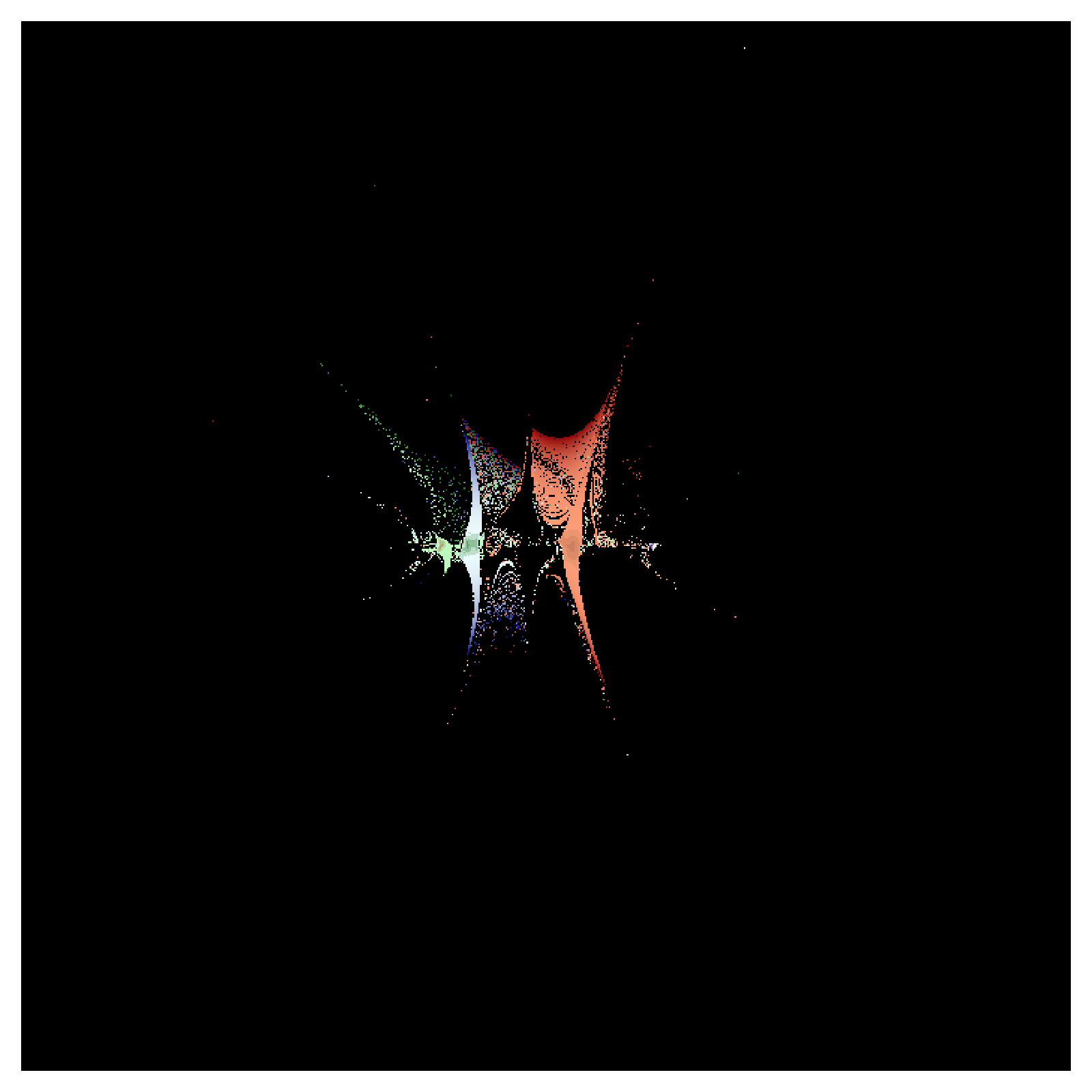"}\\
$f_9$ accelerated  &  \includegraphics[width = 0.2\textwidth]{"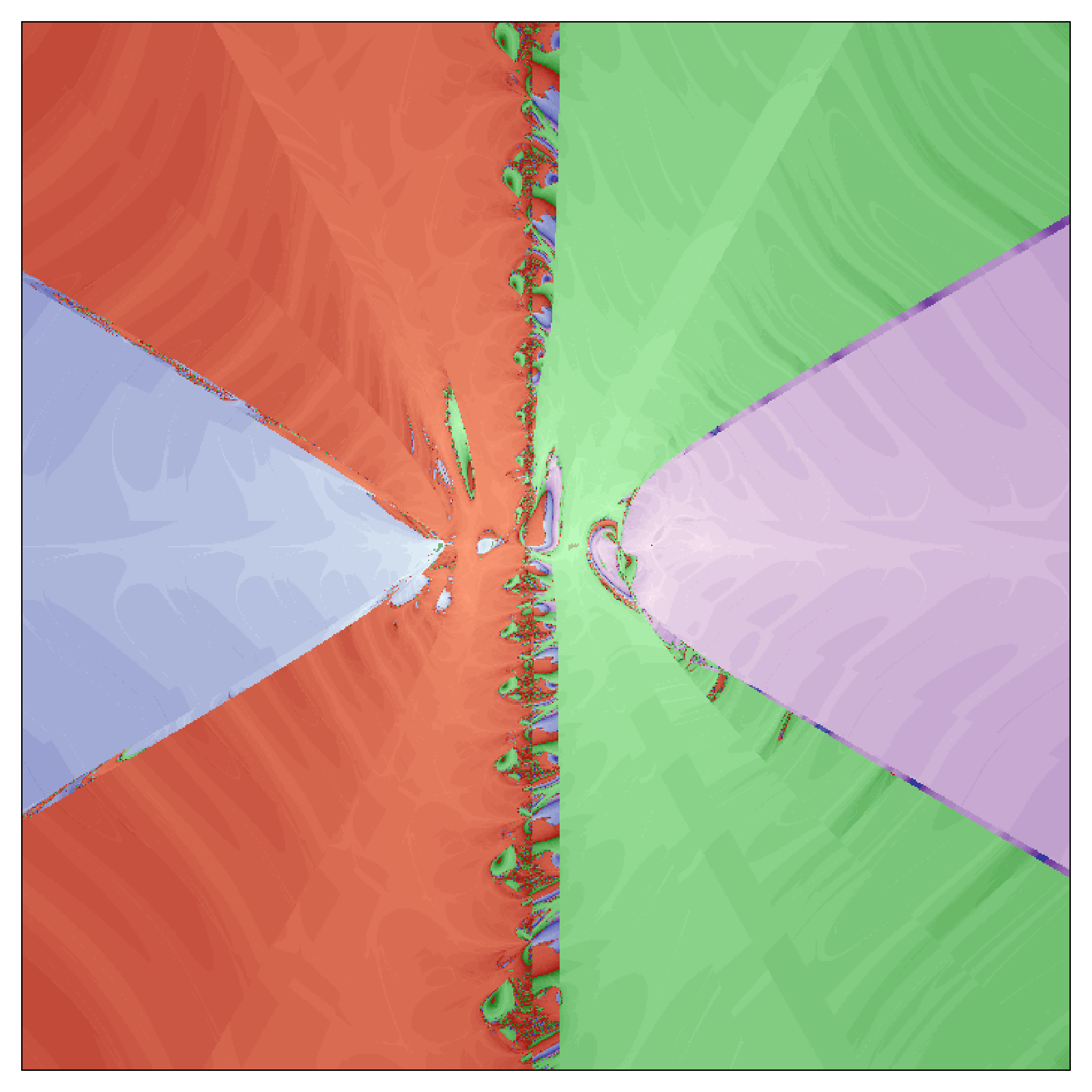"} &
    \includegraphics[width = 0.2\textwidth]{"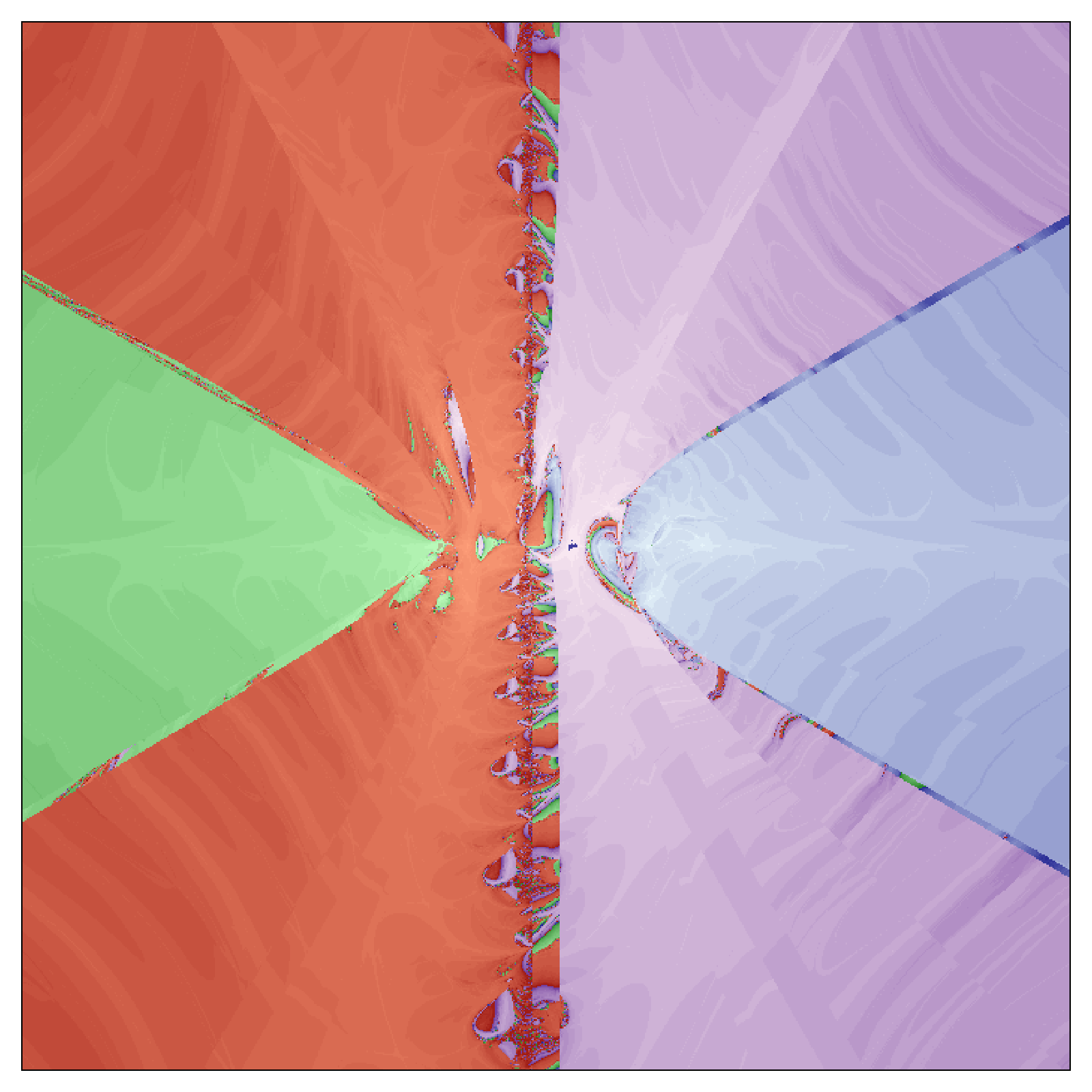"} &
    \includegraphics[width = 0.2\textwidth]{"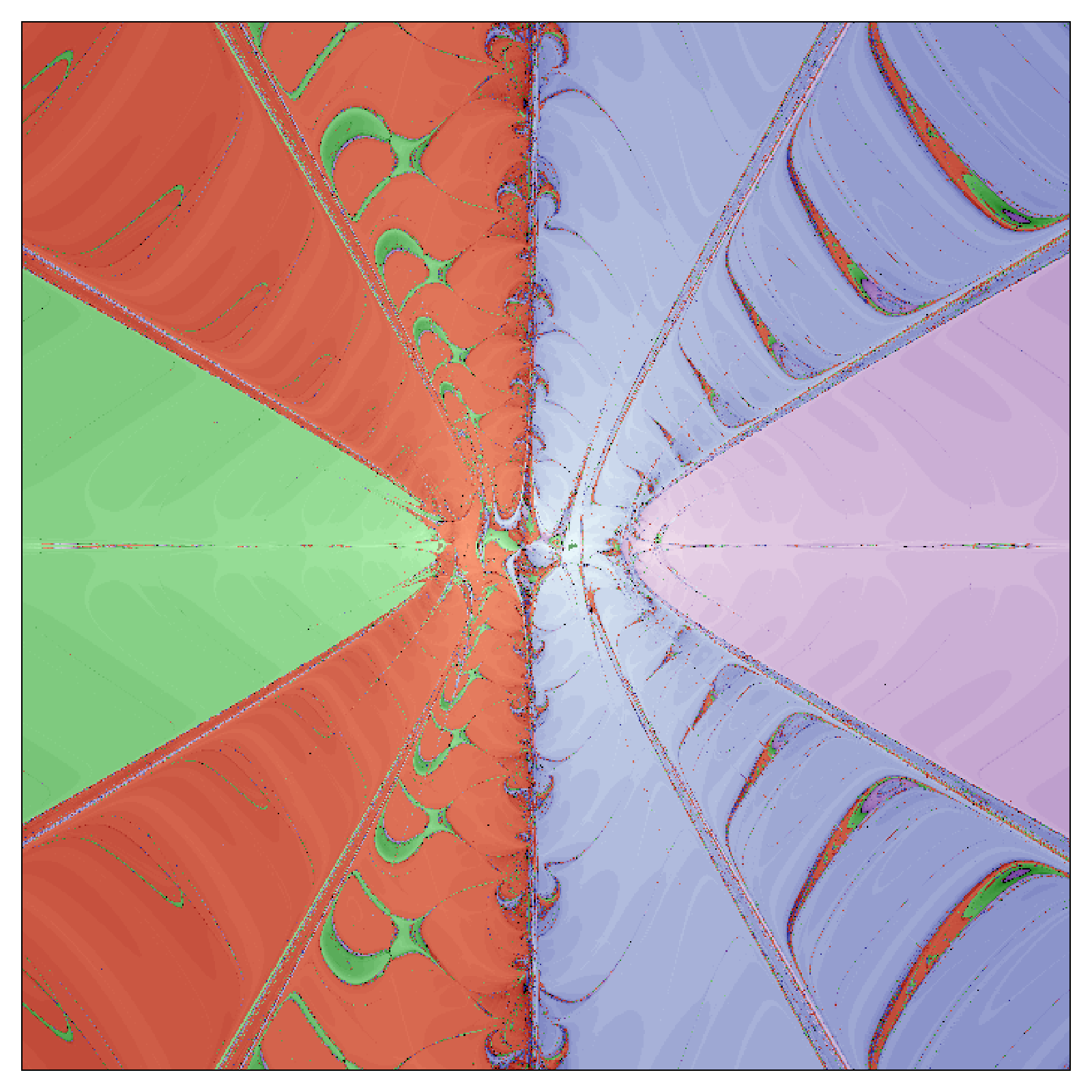"}\\ 
$f_{18}$ normal  &    \includegraphics[width = 0.2\textwidth]{"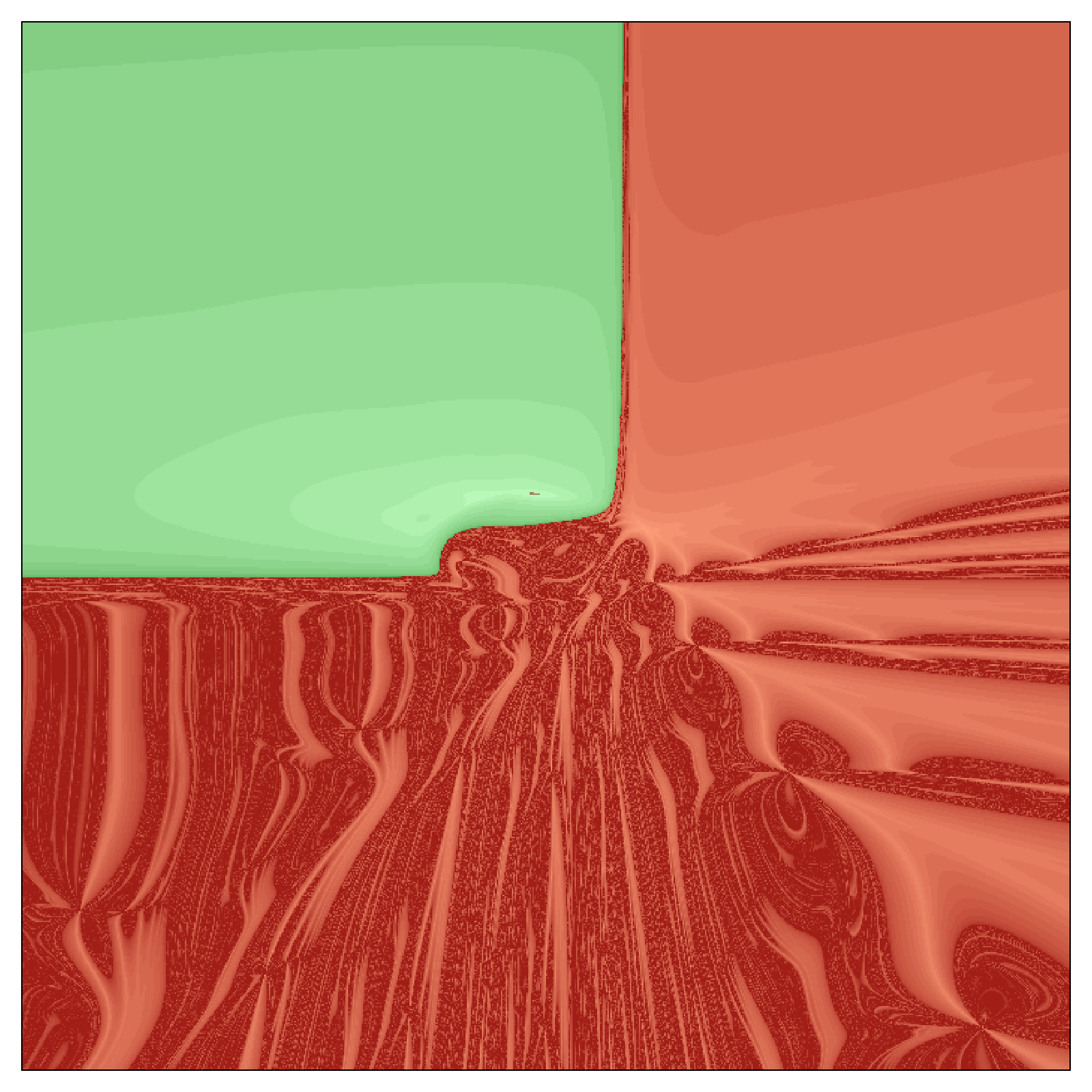"} &
    \includegraphics[width = 0.2\textwidth]{"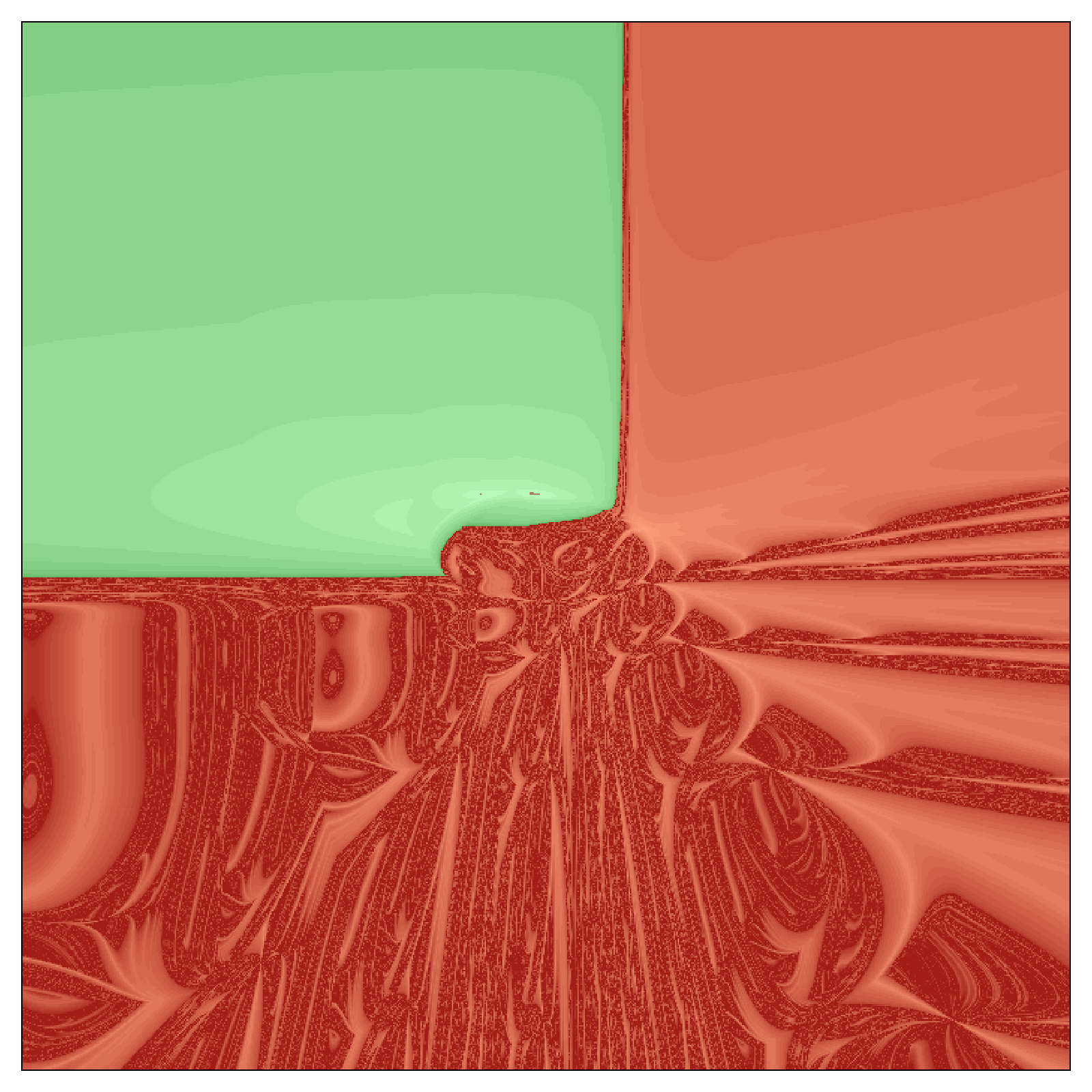"} &
    \includegraphics[width = 0.2\textwidth]{"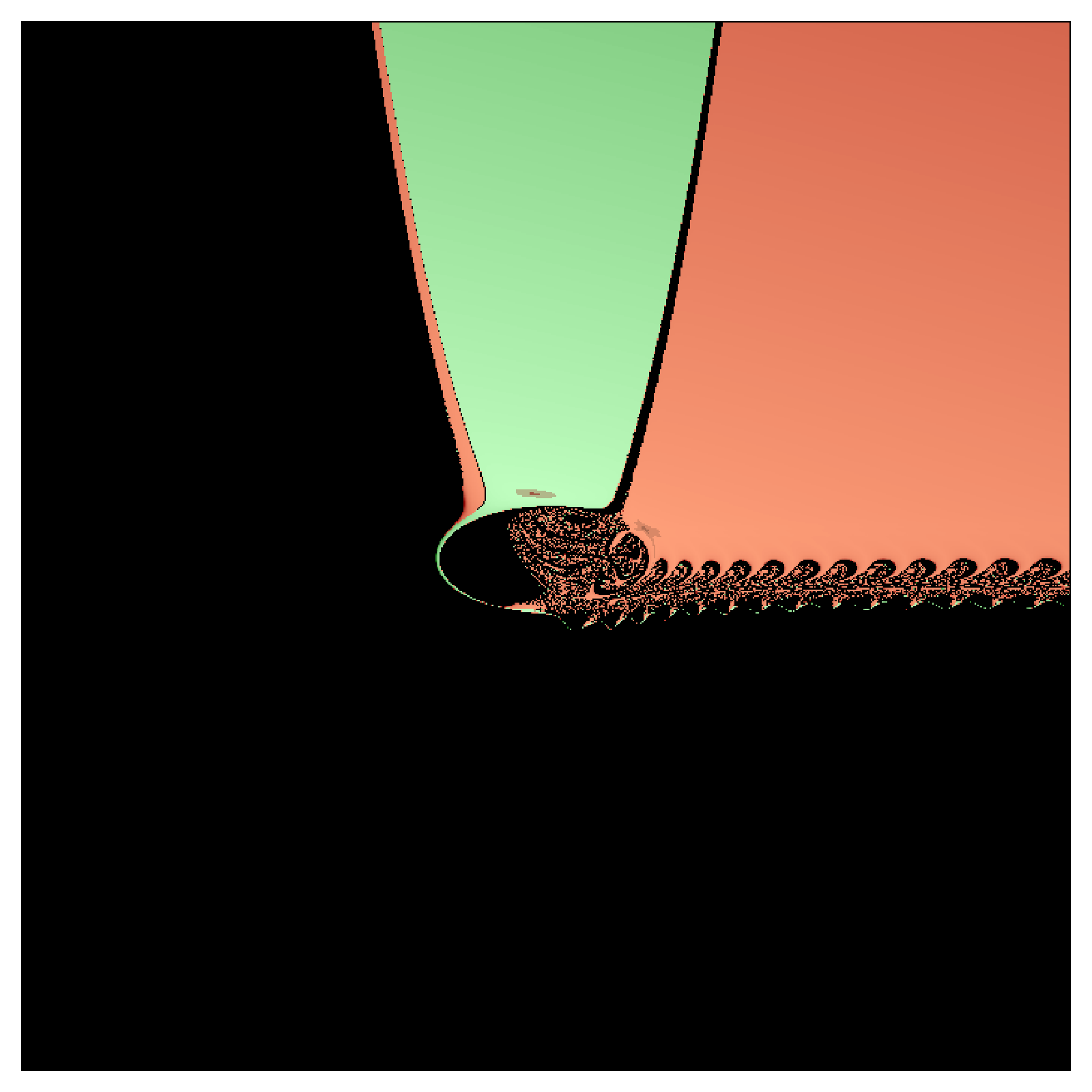"} \\
$f_{18}$ accelerated  &    \includegraphics[width = 0.2\textwidth]{"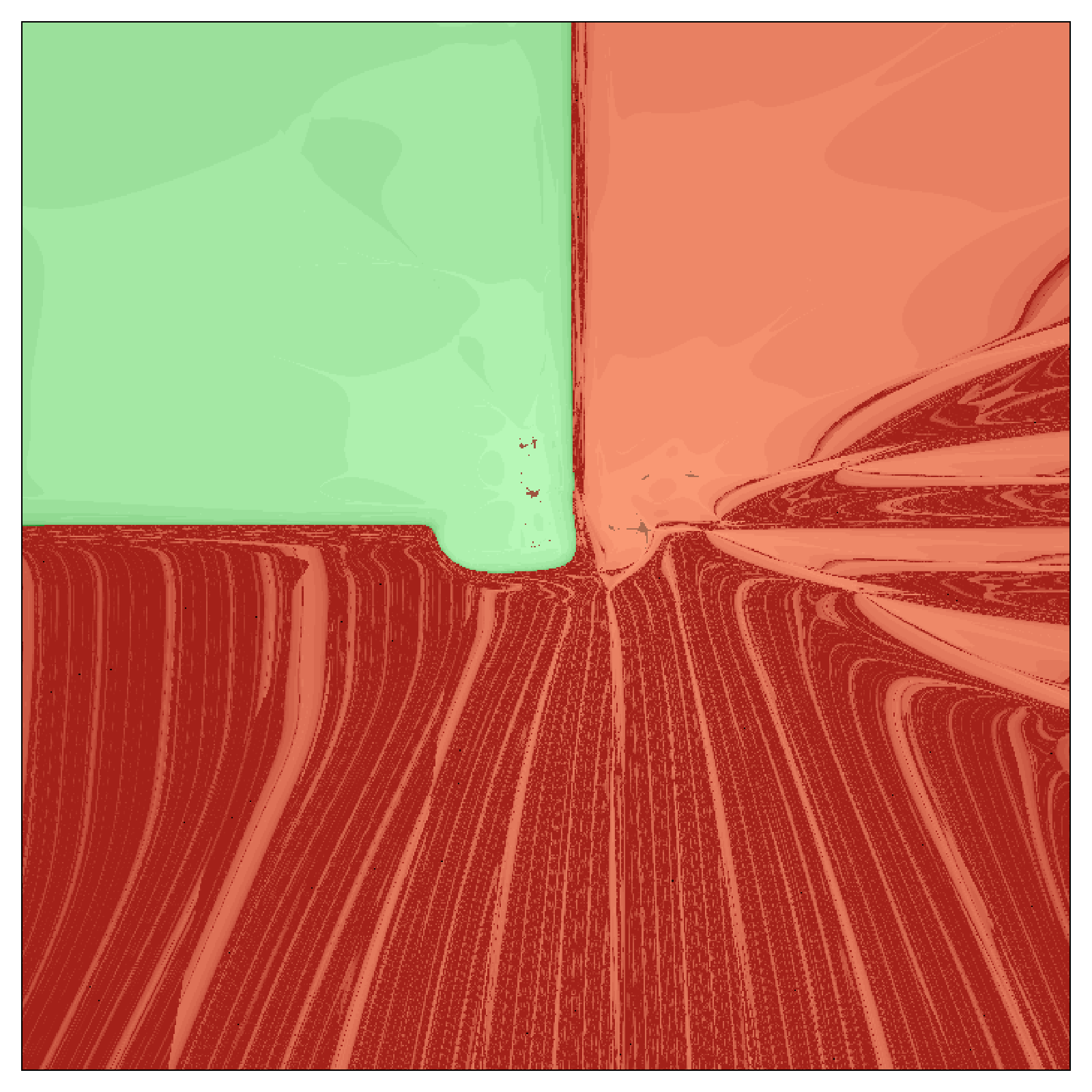"} &
    \includegraphics[width = 0.2\textwidth]{"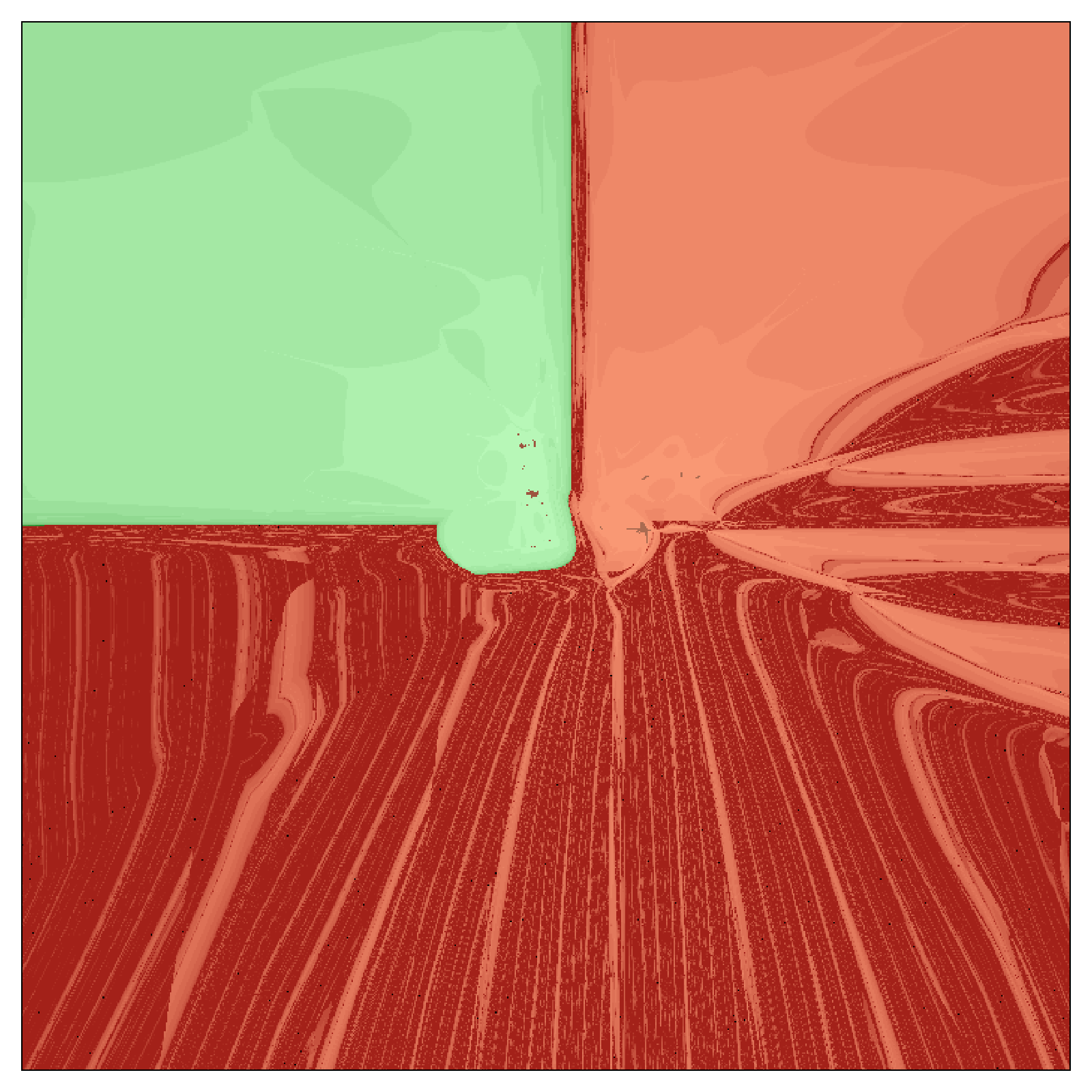"} &
    \includegraphics[width = 0.2\textwidth]{"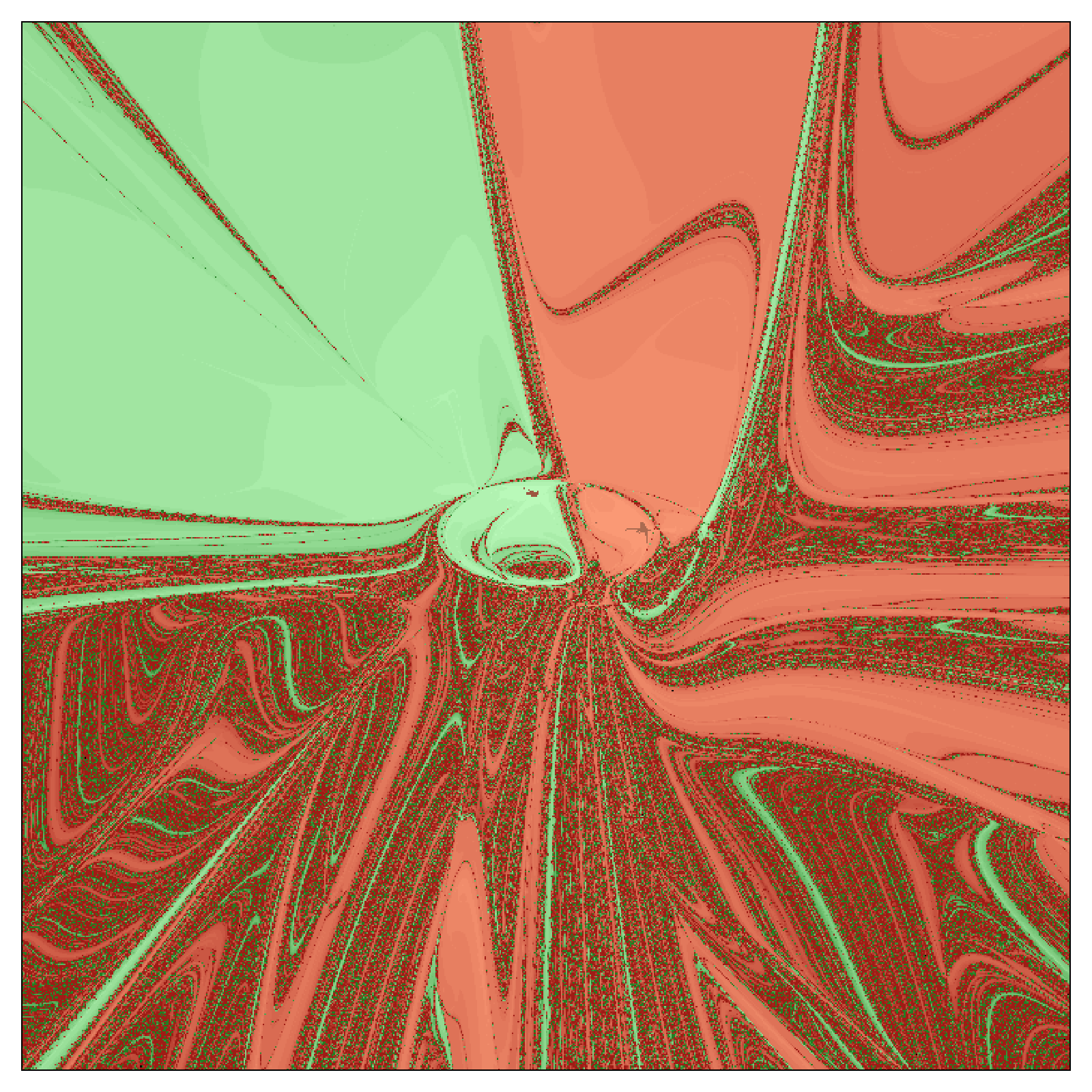"}\\
\end{tblr}
    \caption{\label{fig1} Representation of the basins of attraction of the functions $f_9(z) = (z^2 - 4)(z + 1.5)(z - 0.5)$, and $f_{18}(x,y) =  \left( x^2 - 2x - y + 0.5, x^2 + 4y^2 - 4.0 \right)$. The root-finding algorithm has been computed with the functions $g_i$ described in Sec.~\ref{sec:num}. The classical Steffensen algorithm is simulated in the right column of the panels for the normal and accelerated variations. Plots are computed for a grid of 1000$\times$1000 initial conditions for initial conditions $(x,y)$ in $[-2,2]\times [-2,2]$.}
\end{figure*}

Figure~\ref{fig1} represents the results of the function $f_9$ and $f_{18}$ for all the different methods. Each color in the plot represent a different root except for black representing non-converging orbits. The results for the standard Steffensen method in Fig.~\ref{fig1} show large areas of the basins in black not converging to a root. However, the accelerated methods solve this stability issue. The different algorithms affect the boundary structure in non-trivial ways. The mixing happening in phase space can be exploited as a way to discover different roots more quickly.

\section{Zeros of scalar fields}

\begin{table*}[]
    \centering
    \begin{tabular}{|p{5.5cm}p{1cm} lll | lll | lll |}
    \hline 
 & & \multicolumn{3}{c}{Non conv.} & \multicolumn{3}{|c}{Mean it.}  & \multicolumn{3}{|c|}{T/it.}  \\
 & & $g_1$ & $g_2$ & SM &
    $g_1$ & $g_2$ & SM &
    $g_1$ & $g_2$ & SM\\
 \hline
{\footnotesize $f(x, y) = x^2 + y^2$}
& {\footnotesize (norm.)}
    & \textbf{0.0} & \textbf{0.0} & 51.1 
    & \textbf{12.3} & 12.5 & 52.2 
    & 0.76 & \textbf{0.73} & 1.0 \\
& {\footnotesize (accel.)}
 & 0.0 & 0.0 & 0.0 
    & \textbf{7.4} & \textbf{7.4} & 8.4 
    & 1.54 & \textbf{0.91} & 1.0 \\
 \hline
{\footnotesize $f(x, y) = xy + (x - y) + 100 \sin(y)$}
& {\footnotesize (norm.)}
 & 0.0 & 0.0 & 0.0 
    & \textbf{5.4} & 5.5 & 6.2 
    & \textbf{0.26} & 4.12 & 1.0 \\
& {\footnotesize (accel.)}
 & 0.0 & 0.0 & 0.0 
    & \textbf{4.6} & \textbf{4.6} & 5.1 
    & 2.72 & \textbf{0.88} & 1.0 \\
 \hline
{\footnotesize $f(x, y) = x^2 - (1 - e^y)$}
& {\footnotesize (norm.)}
    & 39.0 & 39.2 & \textbf{12.3} 
    & 14.2 & 14.2 & \textbf{2.4} 
    & 1.7 & 1.62 & \textbf{1.0} \\
& {\footnotesize (accel.)}
    & 0.3 & \textbf{0.2} & 3.9 
    & 30.6 & 30.3 & \textbf{3.3} 
    & 2.72 & \textbf{0.94} & 1.0 \\
 \hline
{\footnotesize $f(x, y, z) = (x + y + z - 3)^2 + (xyz - 1)^2$}
& {\footnotesize (norm.)}
    & 8.9 & \textbf{8.8} & 9.1 
    & 26.4 & 26.3 & \textbf{17.0} 
    & 0.83 & \textbf{0.8} & 1.0 \\
& {\footnotesize (accel.)}
    & 5.5 & 6.1 & \textbf{0.0} 
    & 20.3 & 21.6 & \textbf{15.5} 
    & \textbf{0.93} & 1.08 & 1.0 \\
 \hline
{\footnotesize $f(x, y, z, t) = \sin(x) + 2 \sin(y) + 3 \sin(z) + 4 \cos(t)$}
& {\footnotesize (norm.)}
    & \textbf{91.0} & 91.6 & 91.6 
    & \textbf{5.0} & 4.9 & 4.9 
 & 1.1 & 1.02 & 1.0 \\
& {\footnotesize (accel.)}
    & 78.9 & 77.7 & \textbf{0.0} 
    & 4.2 & \textbf{4.1} & 6.7 
    & 2.22 & 1.26 & \textbf{1.0} \\
 \hline
\end{tabular}

    \caption{\label{tab4}Zeros of scalar fields. The measured metrics are percentage of non-converging initial conditions, average number of iterations and compute time per iteration for $10^4$ sampled initial conditions in the interval $[-10, 10]^d$ where $d$ is the dimension of the function.}
\end{table*}

To complete the possible applications of the modified algorithm, we present an iterative scheme to solve the problem $f(\bar x) = 0$ for a scalar function of $n$ variables. The approximation of gradients to find the optimum of a function is called sub-gradient in the literature. The application of the Newton method to this problem starts with the observation of the first order approximation of $f$: 
\begin{equation}
   f(x_n + \Delta x) \simeq f(x_n) + \Delta x \nabla f(x_n)
\end{equation}
We want an iterate $x_{n+1} = x_n + \Delta x$, such that $f(x_{n+1}) \simeq f(x) + \Delta x \nabla f = 0$. One choice for $\Delta x$ is: 
\begin{equation}
    \Delta x = - \frac{f(x)}{||\nabla f||^2}\nabla f.
\end{equation}
The iterated scheme is then: 
\begin{equation}
x_{n+1} = x_n - \frac{f(x_n)}{||\nabla f||^2}\nabla f.
\end{equation}

If we rewrite the Steffensen iteration using the gradient direction we get:
\begin{equation}
    \bar x_{n+1} = \bar x_n - \alpha_n \nabla \hat f(\bar x_n).
\end{equation}
This is the steepest descent equation but with the estimation of the gradient $\hat f(\bar x_n)$ described in Eq.~\ref{eq:grad_est}. There are several choices for $\alpha_n,$ often called the learning rate or step size. The choice presented above is the Polyak step size~\cite{polyak1969minimization} $\alpha_n = f(\bar x_n)/||\nabla \hat f(\bar x_n)||$. 

Table~\ref{tab4} summarize the results for a few different functions listed. The Steffensen method, in particular the accelerated version, performs well for this task with good results for the functions tested. It is the best candidate for some non-linear functions as shown in the two last rows of the table. The Steffensen method behaves especially well for periodic functions. This is also the case in Tab.~\ref{tab3} for the function $f_4$, $f_5$, $f_7$ and $f_{12}$. The common aspect of these function is the appearance of trigonometric functions. We have no ready answers for this improvement of performance, but we note that the trigonometric functions for these cases are bounded and the derivative for large $x$ is therefore dominated by the other terms in the functions, and is actually better approximated by setting $g(x) = x,$ the identity function. In fact, putting in nontrivial modifying functions like $g_1$ or $g_2$ leads to a worse estimate of the actual derivative than the identity function. Numerically, long orbits appear with very slow convergence to the root, or in the case of the inverse tangent function, the trajectories are unstable for $|x_0| > 1$.

\section{Conclusion}

The literature on root finding methods is enormous, reflecting the critical importance of these methods.  Interesting variants of classical methods are published every month. In this article we point out that simple modifications of a known technique can substantially improve performance with minor side effects. As shown in the results, in many situations the modification is better, especially in terms of the domain of convergence, or at the very least does not degrade the results.  

The extension of the Steffensen algorithm to scalar fields is closely related to optimization problems of finding local minimums without gradient estimations. Possible applications of our technique to this area are an open problem. 

A further improvement of the Steffensen method is to use an annealing approach to speed up the convergence to the root~\cite{jo2024annealing}. This modification requires an additional function evaluation, thus increasing the computational load. While this added complexity can be mitigated with approximations, the possible gains need to be verified empirically. 

\section*{Acknowledgments}
This research was supported [in part] by the Intramural Research Program of the NIH, The National Institute of Diabetes and Digestive and Kidney Diseases (NIDDK) and the European Regional Development Fund (ERDF, EU) under Project No. PID2023-148160NB-I00 (MCIN/AEI/10.13039/501100011033).

\bibliographystyle{elsarticle-num}
\bibliography{steffenson}

\end{document}